\documentclass[12pt]{amsart}
\usepackage{amscd, amssymb, amsfonts}
\usepackage{amsmath,amsthm}
\usepackage{amssymb}
\usepackage{longtable}
\usepackage{amsmath,amscd}
\usepackage[T1]{fontenc}

\usepackage{tikz}
\usetikzlibrary{arrows}
\tikzstyle{block}=[draw opacity=0.7,line width=1.4cm]

  [2002/01/19 v2.2g %
    AMS font definitions%
  ]
\DeclareFontFamily{U}{euf}{}
\DeclareFontShape{U}{euf}{m}{n}{%
  <5><6><7><8><9>gen*eufm%
  <10><10.95><12><14.4><17.28><20.74><24.88>eufm10%
  }{}
\DeclareFontShape{U}{euf}{b}{n}{%
  <5><6><7><8><9>gen*eufb%
  <10><10.95><12><14.4><17.28><20.74><24.88>eufb10%
  }{}

  [2002/01/19 v2.2g %
    AMS font definitions%
  ]
\DeclareFontFamily{U}{msb}{}
\DeclareFontShape{U}{msb}{m}{n}{%
  <5><6><7><8><9>gen*msbm%
  <10><10.95><12><14.4><17.28><20.74><24.88>msbm10%
  }{}

  [2002/01/19 v2.2g %
    AMS font definitions%
  ]
\DeclareFontFamily{U}{msa}{}
\DeclareFontShape{U}{msa}{m}{n}{%
  <5><6><7><8><9>gen*msam%
  <10><10.95><12><14.4><17.28><20.74><24.88>msam10%
  }{}

\newtheorem{theorem}{Theorem}[section]
\newtheorem{lemma}[theorem]{Lemma}

\newtheorem{corollary}[theorem]{Corollary}

\theoremstyle{definition}

\theoremstyle{remark}
\newtheorem{remark}[theorem]{Remark}

\numberwithin{equation}{section}

\begin{document}

\title[]
{On Alzer-Kwong's Identities for Bernoulli polynomials}

\author{Min-Soo Kim}
\address{Division of Mathematics, Science, and Computers, Kyungnam University, 7(Woryeong-dong) kyungnamdaehak-ro, Masanhappo-gu, Changwon-si,
Gyeongsangnam-do 51767, Republic of Korea}
\email{mskim@kyungnam.ac.kr}

\author{Daeyeoul Kim}
\address{Department of Mathematics and Institute of Pure and Applied Mathematics,
Jeonbuk National University, 567 Baekje-daero, Deokjin-gu, Jeonju-si,
Jeollabuk-do, 54896, Republic of Korea.}
\email{kdaeyeoul@jbnu.ac.kr}

\author{Ji Suk So}
\address{Department of Mathematics and Institute of Pure and Applied Mathematics,
	Jeonbuk National University, 567 Baekje-daero, Deokjin-gu, Jeonju-si,
	Jeollabuk-do, 54896, Republic of Korea.}
\email{goleta961@gmail.com}

%\thanks{*Corresponding author}

\begin{abstract}
In this paper, we prove new identities for Bernoulli polynomials that extend Alzer and Kwong's results.
The key idea is to use the Volkenborn integral over $\mathbb Z_p$ of the Bernoulli polynomials
to establish recurrence relations on the integrands. Also, some known identities are obtained by our approach.
\end{abstract}

\subjclass[2000]{11B68, 11S80}
\keywords{Bernoulli polynomials, symmetry properties, Volkenborn integrals}

%\thanks{Received May 24, 2009}

\maketitle

% \linenumbers


\def\ord{\text{ord}_p}
\def\C{\mathbb C_p}
\def\BZ{\mathbb Z}
\def\Z{\mathbb Z_p}
\def\Q{\mathbb Q_p}
\def\wh{\widehat}
\def\ov{\overline}
\def\D{\mathbf D}

\section{Introduction}
\label{Intro}

Bernoulli polynomials play fundamental roles in various branches of mathematics
including combinatorics, number theory, special functions and analysis, see for example \cite{CS,Ko,Sc,Su}.
Let $\mathbb N=\{1,2,\ldots\}$ and $\mathbb N_0=\mathbb N\cup\{0\}.$
The Bernoulli polynomials $B_n(x)$ are usually defined by the generating function
\begin{equation}\label{Eu-pol}
\frac{te^{xt}}{e^t-1}=\sum_{n=0}^\infty B_n(x)\frac{t^n}{n!} \quad  \left( |t|< 2\pi\right).
\end{equation}
And the Bernoulli numbers $B_n$ can be defined by $B_n=B_n(0).$ It is well-known that $B_n=0$ for any odd $n>1.$
These numbers appeared for the first time in Jakob Bernoulli’s book \textit{Ars Conjectandi}, which was published posthumously in 1713.
The polynomials $B_n(x)$ obey the relation $B_n(x)=\sum_{j=0}^n\binom njx^{n-j}B_j.$
Lehmer \cite{Lehm} showed that the Bernoulli polynomials satisfy the relations $B_n(1)=(-1)^nB_n(0)$ and
$B_n(1-x)=(-1)^nB_n(x).$
It is well known that the Bernoulli polynomials have the binomial expansion
$B_n(x+y)=\sum_{j=0}^n\binom njB_j(x)y^{n-j}.$

Alzer and Kwong's \cite{AK} paper was inspired by interesting research note published by Kaneko \cite{Ka} in 1995.
The aim of this paper is to prove the following new identities for Bernoulli polynomials that extend Alzer and Kwong's results
(see \cite{AK}).
In Section \ref{proofs}, we prove the following result.

\begin{theorem}\label{thm-new}
\begin{enumerate}
\item For $m\in\mathbb N$ and $\nu \in \mathbb N_0$ with $0\leq \nu\leq m,$ we have
$$\begin{aligned}
\sum_{\substack{k=0\\ k+m~\text{odd}}}^{m-1}&\binom mk \binom{k+m}{\nu}\binom{k+m-\nu}{m-\nu}B_k(x) \\
&=\frac12\sum_{j=0}^{m-1}(-1)^{j+m+1}\binom m{j+1}\binom{j+m}{\nu}\binom{j+m-\nu}{m-\nu}(j+m+1)x^{j}.
\end{aligned}$$
\item For $m\in\mathbb N$ with $0\leq \nu\leq m-1,$ we have
$$\begin{aligned}
\sum_{\substack{k=0\\ k+m~\text{odd}}}^{m-1}&\binom mk \binom{k+m}{\nu}\binom{k+m-\nu}{m-\nu-1}B_{k+1}(x) \\
&=\frac12\sum_{j=0}^{m}(-1)^{j+m}\binom m{j}\binom{j+m-1}{\nu}\binom{j+m-\nu-1}{m-\nu-1}(j+m)x^{j}.
\end{aligned}$$
\item For $m\in\mathbb N$ with $0\leq \nu\leq m-1$ and $0\leq\ell\leq m-\nu-1,$ we have
$$\begin{aligned}
\sum_{\substack{k=0\\ k+m~\text{odd}}}^{m-1}&\binom mk \binom{k+m}{\nu}\binom{k+m-\nu}{\ell}B_{k+m-\nu-\ell}(x) \\
&=\frac12\sum_{j=0}^{m}(-1)^{j+m}\binom m{j}\binom{j+m-1}{\nu}\binom{j+m-\nu-1}{\ell}\\
&\quad \times (j+m)x^{j+m-\nu-\ell-1}.
\end{aligned}$$
\item For $m\in\mathbb N$ with $0\leq \nu\leq m,$ and $0\leq \ell\leq m-1$  we have
$$\begin{aligned}
\sum_{\substack{k=\ell\\ k+m~\text{odd}}}^{m-1}&\binom mk \binom{k+m}{\nu}\binom{k+m-\nu}{\ell+m-\nu}B_{k-\ell}(x) \\
&=\frac12\sum_{j=\ell}^{m-1}(-1)^{j+m+1}\binom m{j+1}\binom{j+m}{\nu}(j+m+1)x^{j-\ell}.
\end{aligned}$$
\end{enumerate}
\end{theorem}

Setting $\nu=0$ in (3) and (4) of Theorem \ref{thm-new}, we can derive the following corollaries:

\begin{corollary}
\begin{enumerate}
\item For $m\in\mathbb N$ and $0\leq\ell\leq m-\nu-1,$ we have
	$$\begin{aligned}
	\sum_{\substack{k=0\\ k+m~\text{odd}}}^{m-1}&\binom mk \binom{k+m}{\ell}B_{k+m-\ell}(x) \\
	&=\frac12\sum_{j=0}^{m}(-1)^{j+m}\binom m{j}\binom{j+m-1}{\ell}(j+m)x^{j+m-\ell-1}.
	\end{aligned}$$
\item For $m\in\mathbb N$ and $0\leq \ell\leq m-1$  we have
	$$\begin{aligned}
	\sum_{\substack{k=\ell\\ k+m~\text{odd}}}^{m-1}&\binom mk \binom{k+m}{\ell+m}B_{k-\ell}(x) \\
	&=\frac12\sum_{j=\ell}^{m-1}(-1)^{j+m+1}\binom m{j+1}(j+m+1)x^{j-\ell}.
	\end{aligned}$$
\end{enumerate}
\end{corollary}

Setting $\ell=0$ in (3) of Theorem \ref{thm-new}, we can derive the following identity of Alzer and Kwong \cite[Theorem 1]{AK}:

\begin{theorem}[Alzer-Kwong]\label{thm}
Let $m\in\mathbb N$ with $0\leq \nu\leq m.$ Then we have
$$\begin{aligned}
\sum_{\substack{k=0\\ k+m~\text{odd}}}^{m-1}\binom mk&\binom{k+m}{\nu}B_{k+m-\nu}(x) \\
&=\frac12\sum_{j=0}^m(-1)^{j+m}\binom mj\binom{j+m-1}{\nu}(j+m)x^{j+m-\nu-1}.
\end{aligned}$$
\end{theorem}

By setting $x=0$  in (1), (2), (3) and (4) of Theorem \ref{thm-new},
since $0^j=1$ if $j=0$ and $0^j=0$ if $j\in\mathbb N,$
we take the following identities of of Alzer and Kwong \cite[Theorem 2]{AK}:

\begin{theorem}[Alzer-Kwong]\label{thm-1}
\begin{enumerate}
\item For $0\leq \nu\leq m,$ we have
$$\begin{aligned}
\sum_{\substack{k=0\\ k+m~\text{odd}}}^{m-1}\binom mk\binom{k+m}{\nu}\binom{k+m-\nu}{m-\nu}B_k=(-1)^{m+1}\frac{m(m+1)}{2}\binom{m}{\nu}.
\end{aligned}$$
\item For $0\leq \nu\leq m-1,$ we have
$$\sum_{\substack{k=0\\ k+m~\text{odd}}}^{m-1}\binom mk\binom{k+m}{\nu}\binom{k+m-\nu}{m-1-\nu}B_{k+1}=(-1)^m\frac m2\binom{m-1}\nu.$$
\item For $0\leq \ell\leq m-\nu-2,$ we have
$$\sum_{\substack{k=0\\ k+m~\text{odd}}}^{m-1}\binom mk\binom{k+m}{\nu}\binom{k+m-\nu}{\ell}B_{k+m-\nu-\ell}
=0.$$
\item For $0\leq \ell\leq m-1$ and $0\leq \nu\leq m,$ we have
$$\begin{aligned}\sum_{\substack{k=\ell\\ \ell+m~\text{odd}}}^{m-1}\binom mk&\binom{k+m}{\nu}\binom{k+m-\nu}{\ell+m-\nu}B_{k-\ell} \\
&=(-1)^{\ell+m+1}\frac{\ell+m+1}{2}\binom m{\ell+1}\binom{\ell+m}{\nu}.\end{aligned}$$
\end{enumerate}
\end{theorem}

The following identity is due to Wu, Sun and Pan \cite[Theorem 2, (6)]{WSP}.

\begin{theorem}[Wu-Sun-Pan]\label{thm2}
	\begin{align*}
	\sum_{k=0}^{m} \binom mk B_{n+k}(x) = (-1)^{n+m} \sum_{k=0}^{n} \binom nk B_{m+k}(-x),
	\end{align*}
where $m$ and $n$ are positive integers.
\end{theorem}

\begin{theorem} \label{thm3}
We have
	\begin{align*}
	\sum_{j=0}^{m+q} &\binom{m+q}{j}(n+q+j)B_{n+q+j-1}(x) \\
	 &\quad = -(-1)^{m+n} \sum_{k=0}^{n+q}\binom{n+q}{k}(m+q+k)B_{m+q+k-1}(-x),
	\end{align*}
where $q,m$ and $n$ are nonnegative integers and $m+n>0.$
\end{theorem}

\begin{remark}
Substituting $q=1$ and $x=0$ into Theorem \ref{thm3}, we can derive the extension of Kaneko's \cite{Ka} given by
Momiyama \cite{Mom}. It was proved by using a $p$-adic integral over $\mathbb Z_p.$
The Kaneko identity is stated a folllows
$$\sum_{j=0}^{n+1}\binom{n+1}j \tilde{B}_{n+j}=0,$$
where $\tilde{B}_{n}=(n+1)B_n.$
\end{remark}

\begin{theorem} \label{thm4}
We have
	\begin{align*}
	\sum_{j=0}^{m} &\binom{m}{j}\binom{n+j}{\nu}B_{n+j-\nu}(x) \\
	 &\quad = \sum_{k=0}^{n}(-1)^{n-k}\binom{n}{k}\binom{m+k}{\nu}B_{m+k-\nu}(x+1),
	\end{align*}
where $\nu,m$ and $n$ are nonnegative integers and $m+n>0.$
\end{theorem}

Sun \cite[Theorem 1.2, (1.15)]{Su} derived the next identity on Bernoulli polynomials (see also \cite[Theorem 5.1]{CS}).
This identity can be verified by our approach.

\begin{theorem}[Sun] \label{thm5}
We have
	\begin{align*}
	(-1)^m\sum_{j=0}^m\binom mj x^{m-j}B_{n+j}(y)=(-1)^n\sum_{k=0}^n\binom nk x^{n-k}B_{m+k}(z),
	\end{align*}
where $x+y+z=1.$
\end{theorem}

\section{Proofs of Theorem \ref{thm-new}, \ref{thm2} and \ref{thm3}}\label{proofs}

Throughout this section $\Z,\Q$ and $\C$ will, respectively, denote the ring of $p$-adic integers,
the field of $p$-adic rational numbers and the completion of algebraic closure of $\Q.$

For the fundamental properties of $p$-adic integrals and $p$-adic distributions, which are given briefly below,
we may refer the references \cite{Ko,Ro,Sc} and the references cited therein.

The Volkenborn integral of a function $f:\mathbb Z_p\to\C$ is defined by
\begin{equation}\label{-q-e}
\int_{\mathbb Z_p}f(t)dt=\lim_{N\rightarrow\infty}\frac1{p^N}\sum_{j=0}^{p^N-1}f(j)
\end{equation}
and that this limits exists if $f$ is uniformly (or strictly) differentiable on $\Z.$
A function $f:X\to\C$ is uniformly differentiable on $X\subset\C$ (assumed not to have isolated points),
denoted by $f\in UD(\mathbb Z_p),$  if at all points $a\in X$
\begin{equation}\label{st-dif}
\lim_{(x,y)\to(a,a)}\frac{f(x)-f(y)}{x-y}
\end{equation}
exists, the limit being restricted to $x,y\in X,x\neq y$ (see \cite[p. 218]{Ro}).

Let $S\subset\C$ be an arbitrary subset closed under $x\to x+t$ for $t\in\Z$ and $x\in S.$
That is, $S$ could be $\C\backslash\Z,\Q\backslash\Z$ or $\Z.$
Suppose $f:S\to\mathbb C_p$ is strictly differentiable on $S,$ so that for fixed $x\in S$ the function $t\to f(x+t)$ is
uniformly differentiable on $\Z.$

For $f\in UD(\mathbb Z_p),$ the Volkenborn integral
\begin{equation}\label{fer-sim}
F(x)=\int_{\mathbb Z_p}f(x+t)dt, \quad(x\in S)
\end{equation}
is given then satisfied the equation
\begin{equation}\label{fer-sim}
F(x+1)-F(x)=f'(x)
\end{equation}
(see, e.g., \cite[p.~265]{Ro} and \cite[Proposition 55.5(ii)]{Sc}).
From (\ref{fer-sim}), we can be written as
\begin{equation}\label{fer-sim-2}
F(x+q)-F(x+q-1)=f'(x+q-1),
\end{equation}
where $x\in S$ and $q\in \mathbb N.$
%Substitute $q=1,\ldots,n$ and add.
It is easily checked that
%The result is
\begin{equation}\label{fer-sim-3}
F(x+n)-F(x)=\sum_{i=0}^{n-1}f'(x+i)
\end{equation}
with $x\in S$ and $n\in \mathbb N.$

In order to prove Theorem \ref{thm-new}, \ref{thm2} and \ref{thm3}, we need the following lemmas.
It should be noted that the following lemma was obtained by Schikhof in \cite[Proposition 55.7]{Sc}.

\begin{lemma}\label{lem1}
Let $f\in UD(\mathbb Z_p).$  Then we have the functional equation
$$
\int_{\mathbb Z_p}f(-t)dt=\int_{\mathbb Z_p}f(t+1)dt=\int_{\mathbb Z_p}f(t)dt+f'(0).
$$
In particular, if $f$ is an odd function, then
$$\int_{\mathbb Z_p}f(t)dt=-\frac12f'(0).$$
\end{lemma}

It is worth nothing that the Witt's formula for $B_n(x)$ is indeed efficient in deriving recurrence relations
for the ordinary Bernoulli polynomials in an elementary way.

\begin{lemma}[Witt's formula for $B_n(x)$]\label{lem2}
For any $n\in\mathbb N_0,$ we have
$$\int_{\mathbb Z_p}(x+t)^ndt=B_n(x).$$
\end{lemma}

It is known that the Witt's formula for $B_n(x)$ is a power tool in the study of the Bernoulli numbers, Bernoulli polynomials
and its generalization, $p$-adic analytic number theory, etc. (see \cite{Ko,Ro}).

\subsection{Proof of Theorem \ref{thm-new}}

Note that
$$\left( \frac{d}{dt}\right)^{\nu} t^m
=\begin{cases}
\nu!\binom{m}{\nu}t^{m-\nu} &\text{if } m\geq\nu, \\
0 &\text{otherwise}.
\end{cases}$$
For $x\in S,$ we define
$$R(t;x) :=(x+t)^m(x+t-1)^m$$
on $\mathbb Z_p.$
Thus, by the binomial expansion, we find
$$R(t;x)=\sum_{k=0}^m(-1)^{k+m}\binom mk(x+t)^{k+m}$$and
$$R(t+1;x)=\sum_{k=0}^m\binom mk(x+t)^{k+m}.$$
Let $D_t$ be a differentiation operator. Then we set
$$R^{(\nu)}(t;x)= D_{t}^{\nu} R(t;x)= \left( \frac{\partial}{\partial t}\right)^{\nu}R(t;x).$$
It is clear from the definition that
\begin{equation*}
\begin{aligned}
R^{(\nu)}(t+1;x)&-R^{(\nu)}(t;x)\\
&=\nu!\sum_{k=0}^m\binom mk\binom{k+m}\nu(x+t)^{k+m-\nu}\left[1-(-1)^{k+m} \right].
\end{aligned}
\end{equation*}
Since $\left[1-(-1)^{k+m} \right]=0$ if $k$ and $m$ have same parity, we have
\begin{equation}\label{3-1}
\begin{aligned}
R^{(\nu)}(t+1;x)&-R^{(\nu)}(t;x)\\
&=\begin{cases}
	2\nu!\sum_{k=0}^{m-1}\binom mk\binom{k+m}\nu(x+t)^{k+m-\nu} &\text{if $k+m$ odd}, \\
	0 &\text{if $k+m$ even}.
\end{cases}
\end{aligned}
\end{equation}
Similarly, from
\begin{equation} \label{binom}
\binom{j+m}\nu (j+m-\nu)=\binom{j+m-1}\nu(j+m),
\end{equation}
we have
\begin{equation}\label{3-2}
\begin{aligned}
R^{(\nu+1)}(0;x)&=R^{(\nu+1)}(t;x)\biggl|_{t=0} \\
&=\nu!\sum_{j=0}^m(-1)^{j+m}\binom mj\binom{j+m}\nu (j+m-\nu) (x+t)^{j+m-\nu-1} \biggl|_{t=0}\\
&=\nu!\sum_{j=0}^m(-1)^{j+m}\binom mj\binom{j+m-1}\nu(j+m) x^{j+m-\nu-1}.
\end{aligned}
\end{equation}

(1) Let $k+m$ be an odd integer.
Applying $D_{t}^{m-\nu}$ to both side of equation (\ref{3-1}), we see that
\begin{equation}\label{3-1-1}
\begin{aligned}
&D_{t}^{m-\nu}\biggl[R^{(\nu)}(t+1;x)-R^{(\nu)}(t;x)\biggl] \\
&\quad=2\nu!(m-\nu)!\sum_{\substack{k=0\\ k+m~\text{odd}}}^{m-1}\binom mk\binom{k+m}\nu\binom{k+m-\nu}{m-\nu}(x+t)^{k}.
\end{aligned}
\end{equation}
On the other hand, recalling (\ref{3-2}), we get
\begin{equation}\label{3-1-2}
\begin{aligned}
R^{(\nu+1)}(t;x)&=\nu!\sum_{j=1}^{m}(-1)^{j+m}\binom m{j}\binom{j+m}\nu(j+m-\nu)(x+t)^{j+m-\nu-1}\\
&\quad+(-1)^m\nu!\binom m\nu(m-\nu)(x+t)^{m-\nu-1}.
\end{aligned}
\end{equation}
Since $m-\nu>m-\nu-1,$ the second term of right-hand side of (\ref{3-1-2}) gives the identity
$$D_{t}^{m-\nu}\left( (-1)^m\nu!\binom m\nu(m-\nu)(x+t)^{m-\nu-1}\right) =0.$$
Thus, from (\ref{3-1-2}), we have
\begin{equation}\label{3-1-2-1}
\begin{aligned}
D_{t}^{m-\nu}\left[R^{(\nu+1)}(t;x)\right]&=\nu!(m-\nu)!\sum_{j=0}^{m-1}(-1)^{j+m+1}\binom m{j+1}\binom{j+m+1}\nu \\
&\quad\times\binom{j+m-\nu}{m-\nu}(j+m-\nu+1)(x+t)^{j}.
\end{aligned}
\end{equation}
Moreover,
\begin{equation}\label{3-1-2-2}
\begin{aligned}
D_{t}^{m-\nu}\left[R^{(\nu+1)}(0;x)\right]&=\nu!(m-\nu)!\sum_{j=0}^{m-1}(-1)^{j+m+1}\binom m{j+1}\binom{j+m+1}\nu \\
&\quad \times\binom{j+m-\nu}{m-\nu}(j+m-\nu+1)x^{j}.
\end{aligned}
\end{equation}
In particular,
\begin{equation}\label{3-1-3}
\int_{\mathbb Z_p}D_{t}^{m-\nu}\left[R^{(\nu)}(t+1;x)-R^{(\nu)}(t;x)\right]dt=D_{t}^{m-\nu}\left[R^{(\nu+1)}(0;x)\right]
\end{equation}
by using Lemma \ref{lem1}.
On expanding (\ref{3-1-3}) by (\ref{3-1-1}) and (\ref{3-1-2-2}), we obtain
$$\begin{aligned}
\sum_{\substack{k=0\\ k+m~\text{odd}}}^{m-1}&\binom mk\binom{k+m}\nu\binom{k+m-\nu}{m-\nu}
\int_{\mathbb Z_p}(x+t)^{k} dt \\
&=\frac12\sum_{j=0}^{m-1}(-1)^{j+m+1}\binom m{j+1}\binom{j+m}\nu\binom{j+m-\nu}{m-\nu}(j+m+1)x^{j},
\end{aligned}$$
since
$$\binom{j+m+1}\nu(j+m-\nu+1)=\binom{j+m}\nu(j+m+1).$$
Therefore, Part (1) follows from Lemma \ref{lem2}.

(2) From \eqref{3-1} and \eqref{3-2}, we find
\begin{align} \label{3-2-1}
D_{t}^{m-v-1}& \left(R^{(\nu)}(t+1;x)-R^{(\nu)}(t;x)\right) \\
&= 2\nu!(m-v-1)! \sum_{\substack{k=0\\ k+m~\text{odd}}}^{m-1}\binom mk\binom{k+m}\nu\binom{k+m-\nu}{m-\nu-1}(x+t)^{k+1} \nonumber
\end{align}
and
\begin{align} \label{3-2-2}
D_{t}^{m-v-1}\left( R^{(\nu+1)}(t;x)\right)&= \nu!(m-\nu -1)! \sum_{j=0}^{m} (-1)^{j+m} \binom{m}{j} \binom{j+m-1}{\nu}\\
&\quad \times\binom{j+m-\nu-1}{m-\nu-1} (j+m) (x+t)^j. \nonumber
\end{align}
Combining Lemma \ref{lem1} with \eqref{3-2-1} and \eqref{3-2-2}, we get
\begin{align*}
\sum_{\substack{k=0\\ k+m~\text{odd}}}^{m-1} &\binom mk\binom{k+m}\nu\binom{k+m-\nu}{m-\nu-1}\int_{\mathbb Z_p}(x+t)^{k+1} dt  \\
&= \frac{1}{2} \sum_{j=0}^{m} (-1)^{j+m}  \binom{m}{j} \binom{j+m}{\nu} \binom{j+m-1}{\nu}(j+m) x^j.
\end{align*}
So we have Part (2) by Lemma \ref{lem2}.

(3) From \eqref{3-1} and \eqref{3-2}, we known that
\begin{align} \label{3-3-1}
D_{t}^{\ell} (R^{(\nu)}&(t+1;x)-R^{(\nu)}(t;x)) \\
&= 2\nu ! \ell! \sum_{\substack{k=0\\ k+m~\text{odd}}}^{m-1}\binom mk\binom{k+m}\nu\binom{k+m-\nu}{\ell}(x+t)^{k+m-\nu -\ell} \nonumber
\end{align}
and
\begin{equation}
\begin{aligned} \label{3-3-2}
D_{t}^{\ell}\left( R^{(\nu+1)}(t;x)\right)&= \nu! \ell! \sum_{j=0}^{m} (-1)^{j+m} \binom{m}{j} \binom{j+m}{\nu}\binom{j+m-\nu-1}{\ell}  \\
&\quad \times(j+m-\nu) (x+t)^{j+m-\nu-\ell-1}.
\end{aligned}
\end{equation}
Combining Lemma \ref{lem1} with \eqref{3-3-1}, \eqref{3-3-2}
and using \eqref{binom}, we get
\begin{align*}
&\sum_{\substack{k=0\\ k+m~\text{odd}}}^{m-1} \binom mk\binom{k+m}\nu\binom{k+m-\nu}{\ell}\int_{\mathbb Z_p}(x+t)^{k+m-\nu-\ell} dt  \\
&= \frac{1}{2} \sum_{j=0}^{m} (-1)^{j+m}  \binom{m}{j}\binom{j+m-1}{\nu}(j+m) \binom{j+m-\nu-1}{\ell} x^{j+m-\nu-\ell-1},
\end{align*}
which implies Part (3) by using Lemma \ref{lem2}.

(4) By \eqref{3-1} and \eqref{3-1-2}, we have
\begin{align} \label{3-4-1}
&D_{t}^{\ell+m-\nu} \left(R^{(\nu)}(t+1;x)-R^{(\nu)}(t;x)\right) \\
&\quad = 2\nu ! (\ell+m-\nu)! \sum_{\substack{k=0\\ k+m~\text{odd}}}^{m-1}\binom mk\binom{k+m}\nu\binom{k+m-\nu}{\ell+m-\nu}(x+t)^{k-\ell} \nonumber
\end{align}
and
\begin{align} \label{3-4-2}
D_{t}^{\ell+m-\nu}\left( R^{(\nu+1)}(t;x)\right)&= \nu !(\ell+m-\nu)! \sum_{j=\ell}^{m-1} (-1)^{j+m+1} \binom{m}{j+1} \binom{j+m+1}{\nu} \\
&\quad \times(m+j-\nu+1) (x+t)^{j-\ell}. \nonumber
\end{align}
Combining Lemma \ref{lem1} with \eqref{3-4-1} and \eqref{3-4-2}, we get
\begin{align*}
&\sum_{\substack{k=0\\ k+m~\text{odd}}}^{m-1} \binom mk\binom{k+m}\nu\binom{k+m-\nu}{\ell+m-\nu}\int_{\mathbb Z_p}(x+t)^{k-\ell} dt  \\
&\quad= \frac{1}{2} \sum_{j=\ell}^{m-1} (-1)^{j+m+1}  \binom{m}{j+1}\binom{j+m}{\nu}(j+m+1)  x^{j-\ell},
\end{align*}
the right hand side following from
$$\binom{j+m+1}{\nu} (m+j-\nu+1) = \binom{j+m}{\nu}(j+m+1).$$
Hence, Part (4) is obtained by Lemma \ref{lem2}.

\subsection{Proof of Theorem \ref{thm2}}

Let $m$ and $n$ be positive integers. For $x\in S,$ we define
$$G(t;x):=(-1)^m(x+t)^m(x+t-1)^n$$
on $\mathbb Z_p.$
By the binomial expansion, the formula $G(-t;x)$ and $G(t+1;x)$ can be rewritten as
$$G(-t;x)= (-1)^{n}\sum_{k=0}^n\binom nk (-x+t)^{m+k}$$and
$$G(t+1;x)=(-1)^{m}\sum_{k=0}^m\binom mk(x+t)^{n+k}.$$
Applying Lemma \ref{lem1} to the above two equations, we have
\begin{align*}
(-1)^m\sum_{k=0}^{m} \binom mk \int_{\mathbb Z_p}(x+t)^{n+k} dt=(-1)^n\sum_{k=0}^{n} \binom mk \int_{\mathbb Z_p}(-x+t)^{m+k} dt.
\end{align*}
The proof now follows directly from Lemma \ref{lem2}.

\subsection{Proof of Theorem \ref{thm3}}

Let $m,n$ and $q$ be nonnegative integers with $m+n>0.$
For $x\in S,$ we define
$$H(t;x) :=(x+t)^{m+q}(x+t-1)^{n+q}+(-1)^{m+n}(-x+t)^{n+q}(-x+t-1)^{m+q}$$
on $\mathbb Z_p.$
By the binomial expansion, the formula $H(-t;x)$ and $H(t+1;x)$ can be rewritten as
$$H(-t;x)= (-1)^{m+n}(-x+t)^{m+q}(-x+t+1)^{n+q}+(x+t)^{n+q}(x+t+1)^{m+q}$$and
$$H(t+1;x)=(x+t+1)^{m+q}(x+t)^{n+q}+(-1)^{m+n}(-x+t+1)^{n+q}(-x+t)^{m+q}.$$

Since $H(-t;x)=H(t+1;x)$, we have
$$-H'(-t;x) = H'(t+1;x),$$
where $H'=\frac{\partial H}{\partial t}.$
Hence, by Lemma \ref{lem1}, we obtain
\begin{align*}
-\int_{\mathbb Z_p} H'(t+1;x) dt = \int_{\mathbb Z_p} H'(-t;x) dt = \int_{\mathbb Z_p} H'(t+1;x) dt.
\end{align*}
By the above equation, we have
\begin{align} \label{H'}
\int_{\mathbb Z_p} H'(t+1;x) dt = 0.
\end{align}
In particular,
\begin{equation}\label{H'-1}
\begin{aligned}
H'(t+1;x) = &\sum_{j=0}^{m+q} \binom{m+q}{j}(n+q+j)(x+t)^{n+q+j-1} \\
&\quad +(-1)^{m+n} \sum_{k=0}^{n+q}\binom{n+q}{k}(m+q+k)(-x+t)^{m+q+k-1}.
\end{aligned}
\end{equation}
From equation \eqref{H'} and \eqref{H'-1}, we get
\begin{align*}
0&=\int_{\mathbb Z_p} H'(t+1;x) dt \\
 &=\sum_{j=0}^{m+q} \binom{m+q}{j}(n+q+j)\int_{\mathbb Z_p}(x+t)^{n+q+j-1}dt \\
&\quad +(-1)^{m+n} \sum_{k=0}^{n+q}\binom{n+q}{k}(m+q+k)\int_{\mathbb Z_p}(-x+t)^{m+q+k-1} dt.
\end{align*}
Thus the result follows from Lemma \ref{lem2}.

\subsection{Proof of Theorem \ref{thm4}}

We begin by proving some binomial coefficient identities.
We start with the identity
$$
\begin{aligned}
0&=(x+t)^n(x+t+1)^m-(x+t+1)^m(x+t+1-1)^n \\
&=\sum_{j=0}^m\binom mj(x+t)^{n+j}-\sum_{k=0}^n\binom nk(-1)^{n-k}(x+t+1)^{m+k},
\end{aligned}
$$
or equivalently,
\begin{equation}\label{thm4-p1}
\sum_{j=0}^m\binom mj(x+t)^{n+j}=\sum_{k=0}^n\binom nk(-1)^{n-k}(x+t+1)^{m+k}.
\end{equation}
Differentiating both sides of the identity (\ref{thm4-p1}) with respect to $t,$ $\nu$ times,
leads to following relation
\begin{equation}\label{thm4-p2}
\sum_{j=0}^m\binom mj\binom{n+j}{\nu}(x+t)^{n+j-\nu}=\sum_{k=0}^n(-1)^{n-k}\binom nk\binom{m+k}{\nu}(x+t+1)^{m+k-\nu}.
\end{equation}
Applying Lemma \ref{lem1} to (\ref{thm4-p2}), we have
$$
\begin{aligned}
\sum_{j=0}^m\binom mj&\binom{n+j}{\nu}\int_{\mathbb Z_p}(x+t)^{n+j-\nu}dt \\
&=\sum_{k=0}^n(-1)^{n-k}\binom nk\binom{m+k}{\nu}\int_{\mathbb Z_p}(x+t+1)^{m+k-\nu}dt.
\end{aligned}
$$
From Lemma \ref{lem2}, this gives the formula of our theorem.

\subsection{Proof of Theorem \ref{thm5}}

Let $x+y+z=1.$
For $x\in S,$ we define
$$L(t;x) :=(-1)^m(y-1+t)^n(y+x-1+t)^m$$
on $\mathbb Z_p.$
By the binomial expansion, we have
\begin{equation}\label{thm5-p1}
L(t+1;x)=(-1)^m\sum_{j=0}^m \binom mj x^{m-j}(t+y)^{n+j},
\end{equation}
and since $x+y+z=1,$
\begin{equation}\label{thm5-p2}
\begin{aligned}
L(-t;x)&=(-1)^n(t+x+z)^n(t+z)^m \\
&=(-1)^n\sum_{k=0}^n \binom nk x^{n-k}(t+z)^{m+k}.
\end{aligned}
\end{equation}
Applying Lemma \ref{lem1} to (\ref{thm5-p1}) and (\ref{thm5-p2}), we have
\begin{equation}\label{thm5-p3}
\begin{aligned}
(-1)^m\sum_{j=0}^m \binom mj& x^{m-j}\int_{\mathbb Z_p}(y+t)^{n+j}dt \\
&=(-1)^n\sum_{k=0}^n \binom nk x^{n-k}\int_{\mathbb Z_p}(z+t)^{m+k}dt.
\end{aligned}
\end{equation}
Therefore, the result follows from Lemma \ref{lem2}.

\bibliography{central}

\begin{thebibliography}{00}


\bibitem{AK} H. Alzer and M.K. Kwong,
\textit{Identities for Bernoulli polynomials and Bernoulli numbers},
Arch. Math. (Basel)  \textbf{102}  (2014), no. 6, 521--529.

\bibitem{CS} W.Y.C. Chen and L.H. Sun,
\textit{Extended Zeilberger's algorithm for identities on Bernoulli and Euler polynomials},
J. Number Theory \textbf{129} (2009), 2111--2132.

\bibitem{Ka} M. Kaneko,
\textit{A recurrence formula for the Bernoulli numbers},
Proc. Japan Acad. Ser. A. Math. Sci. \textbf{71} (1995), 192--193.

\bibitem{Ko} N. Koblitz,
\textit{$p$-adic Numbers, $p$-adic Analysis, and Zeta-Functions},
New York: Springer-Verlag, 1984.

\bibitem{Lehm} D.H. Lehmer,
\textit{A new approach to Bernoulli polynomials},
Amer. Math. Monthly \textbf{95} (1988) 905--911.

\bibitem{Mom} H. Momiyama,
\textit{A new recurrence formula for Bernoulli numbers},
Fibonacci Quart.  \textbf{39}  (2001), no. 3, 285--288.

\bibitem{Ro} A.M. Robert,
\textit{A course in $p$-adic analysis},
Graduate Texts in Mathematics, 198, Springer-Verlag, New York, 2000.

\bibitem{Sc} W.H. Schikhof,
\textit{Ultrametric calculus. An introduction to $p$-adic analysis},
Cambridge Studies in Advanced Mathematics, \textbf{4}, Cambridge University Press, Cambridge, 1984.

\bibitem{Su} Z.-W. Sun,
\textit{Combinatorial identities in dual sequences},
European J. Combin. \textbf{24} (2003), no. 6, 709--718.

\bibitem{WSP}K.-J. Wu, Z.-W. Sun and H. Pan,
\textit{Some identities for Bernoulli and Euler polynomials},
Fibonacci Quart. \textbf{42} (2004) 295--299.


\end{thebibliography}

\end{document}